\documentclass[reqno,11pt]{amsart}

\usepackage[sort&compress,numbers]{natbib}
\usepackage{array,longtable} 
\usepackage{tabularray}
\usepackage{amsfonts}
\usepackage{mathdots,amssymb,graphicx,amssymb,amsmath,amsopn, mathtools,mathabx}
\usepackage[utf8]{inputenc}
\usepackage{enumitem}
\usepackage{float}
\usepackage{tikz-cd}
\usepackage{scalerel}
\usepackage{pgfplots}
\usepackage{amscd}
\usepackage{latexsym}
\usepackage{xcolor}
\usepackage{amscd}
\usepackage{graphicx}
\usepackage{subcaption}
\usepackage{pb-diagram}
\usepackage{mathrsfs}
\usepackage{accents}
\usepackage[textwidth=15cm,textheight=22.5cm,headheight=0.6cm,headsep=1cm,centering]{geometry}
\usepackage{fancyhdr}
\usepackage{float}
\usepackage[bookmarks=true,colorlinks=true,linkcolor=airforceblue,citecolor=brightmaroon,plainpages=tfalse,pdftex]{hyperref}
\hypersetup{bookmarksopen=true,pdfview=fitbh}

\definecolor{brightmaroon}{rgb}{0.76, 0.13, 0.28}
\definecolor{airforceblue}{rgb}{0, 0.25, 0.77}
\definecolor{myOrange}{rgb}{1,0.5,0}
\definecolor{brightmaroon}{rgb}{0.76, 0.13, 0.28}
\definecolor{airforceblue}{rgb}{0, 0.4, 0.66}

\allowdisplaybreaks
\pgfplotsset{compat=1.18}

\theoremstyle{plain}

\newtheorem{teo}{Theorem}[section]

\newtheorem{pro}[teo]{Proposition}
\newtheorem{coro}[teo]{Corollary}
\newtheorem{defi}[teo]{Definition}

\numberwithin{equation}{section}

\newcommand{\D}{\mathbf{D}}
\newcommand{\J}{\mathbf{J}}
\newcommand{\Px}{\mathbb{P}}

\newcommand{\un}{\mathbf{u}}

\newcommand{\wtP}{\widetilde{P}}

\numberwithin{equation}{section}
\newcommand{\prodint}[1]{\left\langle{#1}\right\rangle}

\newcommand{\C}{\mathbb{C}}
\usepackage{color}

\newcommand*\elln{\ensuremath{\boldsymbol\ell}}

\newcommand*\pFq[2]{{}_{#1}F_{#2}}

\begin{document}
	
	\title[A note on Laguerre truncated polynomials and quadrature formula]{A note on Laguerre truncated polynomials and quadrature formula}

	\author[J. C García-Ardila, F. Marcellán]
	{Juan C. García-Ardila,  Francisco Marcellán}

	\address[J. C. Garc\'ia-Ardila]{Departamento de Matem\'atica Aplicada a la Ingenier\'ia Industrial \\Universidad Polit\'ecnica de Madrid\\ Calle Jos\'e Gutierrez Abascal 2, 28006 Madrid, Spain.}
	\email{juancarlos.garciaa@upm.es}

	\address[F. Marcellán]{Departamento de Matemáticas, Universidad Carlos III de Madrid, Leganés, Spain} \email{pacomarc@ing.uc3m.es}
	
	\thanks{ }
	
	\date{\today}
	\begin{abstract}
In this contribution we deal with Gaussian quadrature rules based on orthogonal polynomials associated with a weight function $w(x)= x^{\alpha} e^{-x}$ supported on an interval $(0,z)$, $z>0.$ The modified Chebyshev algorithm is used in order to test the accuracy in the computation of the coefficients of the three-term recurrence relation, the zeros and weights, as well as the dependence on the parameter $z.$
		
	\end{abstract}
	
	\maketitle

\begin{center}
Dedicated to Academician Professor Gradimir Milovanovi\'c

on the occasion of his 75th birthday.
\end{center}
\section {Introduction}
Integration with respect to a probability measure supported on the real line  is a topic that belongs to analysis, and so is the evaluation or approximation of integrals. According to the seminal work by  Gauss, one is led to orthogonal polynomials with respect to such a measure, It is well known that the zeros of such polynomials are real, simple and located in the interior of the convex hull of the support of the measure. Another important fact is that the zeros of two consecutive polynomials interlace, i. e., between two zeros of the orthogonal polynomial of degree $n+1$ you have exactly one zero of the polynomial of degree $n$ (see \cite{Chihara}, among others). \\

The key elements  in the study of quadrature formulas are the nodes and weights in order to have the optimal degree of exactness. In fact, it is very well known that the Gaussian quadrature formulas provide the optimal degree of exactness for polynomials of degree at most $2n-1$ when you deal with $n$ nodes. Indeed, they are the zeros of the polynomial of degree $n$ associated with the measure appearing in the integral.\\

The connection between Gaussian quadrature rules and the algebraic eigenvalue problem was, if not discovered, then certainly exploited in the now classical and widely cited paper \cite{Go69}. Indeed,  the connection between the zeros of orthogonal polynomials and the eigenvalues of the truncated Jacobi matrices representing the symmetric operator of multiplication by $x$ for the standard inner product defined by the probability measure is the key point.  Thus, we have an interplay between classical analysis and numerical linear algebra as pointed out in \cite{Ga02}. For \textit{Gauss-Radau} and \textit{Gauss-Lobatto quadrature formulas,} the nodes are related with zeros of orthogonal polynomials  with respect to polynomial perturbations of the probability measure of degrees one and two, respectively. Those zeros are eigenvalues of Jacobi matrices that are defined from Darboux transformations of the Jacobi matrix associated with the initial probability measure. On the other hand, in~\cite{Ga82}, the entries of Jacobi matrices associated  with generalized Christoffel transformations (perturbations of probability measures by positive rational functions) and their numerical implementation have been analyzed, where the connection with the QD-algorithm is  also stated.\\

In the literature, Gaussian quadrature formulas  for classical orthogonal polynomials (Hermite, Laguerre and Jacobi) and efficient algorithms to estimate the reminders in the approximation of integral of continuous and analytic functions have been implemented (see \cite{G83} among others).\\

For nonclassical measures, Gaussian quadrature formulas have been  exhaustively  ana\-ly\-zed in the literature  from the pioneering work  \cite{D76} for an absolutely continuous measure $ e^{-tx^2} dx$ supported in the interval $(-1,1)$ where $t$ is a positive real number.  They constitute a basic tool in  the computation of  integrals related to electron repulsion in molecular quantum mechanics. The corresponding sequences of orthogonal polynomials are known in the literature as Rys polynomials. We must point out that such a perturbation of the Lebesgue measure supported in $(-1,1)$ yields to  an example of a Volterra integrable system.  On the other hand, by using a change of variables, the orthogonality with respect to a truncated Normal distribution supported in  a symmetric interval of the real line appears in a natural way. The corresponding  sequences of  orthogonal polynomials  appear in the study of Gaussian unitary ensembles with two discontinuities. Analytic properties of such orthogonal polynomials (called truncated Hermite polynomials) have been studied in \cite{DM22}.  Other examples of  Gauss quadrature formulas for nonclassical measures appear in \cite{M15, MV22}.\\

More recently,   in  \cite{MV22}  the authors have considered the problem of construction of quadrature formulas of Gaussian type on $[-1, 1],$ concerning the following two-parametric weight function
$w(x; z) = e^{-zx^2}(1-x^2)^{\lambda -1/2}$, $\lambda >-1/2.$
This weight function is a generalization of the \textit{Rys quadrature}, and it can be obtained by transforming quadratures on $[-1,1]$ with~$N$ nodes to quadratures on $(0,1)$ with only $(N+1)/2$ nodes. Such an approach provides a stable and very efficient numerical construction. Notice that a Volterra lattice associated with the symmetric Jacobi (Gegenbauer) measure appears in a natural way. It can be also considered  as a continuous perturbation of the Gaussian unitary ensemble with two jumps. \\

In \cite{DGM23}, the authors deal with sequences of  polynomials orthogonal with respect to the truncated Gamma distribution, i. e., the absolutely continuous measure $x^{\alpha} e^{-x} dx$ su\-ppor\-ted on an interval $(0,z),$ $z>0.$ Analytic properties of such polynomials (known as truncated Laguerre polynomials)  as well as of their zeros have been studied therein. Notice that a symmetrization problem yields the truncated generalized Hermite polynomials whcih have been also analyzed in \cite{DGM23}. The aim of the present contribution is to consider Gaussian quadrature formulas for the truncated Gamma distribution as well as to provide an efficient algorithm for the accuracy of them in terms of the parameter~$z$. We emphasize our work on the analysis of  computational methods to deal with the nodes and weights of such Gaussian quadrature formulas.\\

The structure of the manuscript is the following. In Section \ref{S2} we summarize the basic background about  Gaussian quadrature formulas and we explain the modified Chebyshev algorithm based on generalized moments of the measure defined by a new sequence of orthogonal polynomials in order to improve the numerical stability of the eigenvalue problem of the Jacobi matrix. In Section \ref{S3}  we study the generalized moments of the truncated Gamma distribution taking into account shifted Jacobi polynomials as supporting sequence. Their expressions  are deduced. Finally, in Section \ref{S4}, some numerical tests are presented in order to get the coefficients of the three term recurrence relation that the truncated Laguerre polynomials satisfy and their dependence of $z.$ The nodes and weights in the corresponding Gaussian quadrature formulas are described in an illustrative case for choices of $\alpha$ and $z.$
\section{Basic background}\label{S2}
Let $\Px$ be the set of polynomials with complex coefficients. We will denote by $\Px_{n}$ the linear subspace of polynomials of degree at most $n.$ The algebraic dual of $\Px$, denoted by $\Px^*$, is defined as the set of all linear mappings from $\Px$ into $\C$, that is,
$$
\mathbb{P}^* \,=\,\{\mathbf{u}:\mathbb{P}\rightarrow \mathbb{C}\,:\,\mathbf{u}\text{ is linear} \}.
$$
The image of a polynomial $p$ under $\un$ will be expressed using duality brackets as $\prodint{\un,p}$. The elements of $\Px^*$ are usually called \textit{linear functionals}.
Any linear functional $\un$ is completely defined by the values $\un_n=:\prodint{\un,x^n}$, $n\ge 0,$  where  $\un_n$ is called  the $n$-th \textit{moment} of the linear functional $\mathbf{u}$. It is extended by linearity to all polynomials.\\

We turn our attention to orthogonality with respect to a linear functional.
The linear functional $\un$ is said to be \textit{quasi-definite} (resp. \textit{positive-definite}) if every leading principal submatrix of the Hankel matrix $H=(\un_{i+j})_{i , j=0}^{\infty}$ is nonsingular (resp. positive-definite).  In such a situation, there exists a sequence of monic polynomials $(P_n)_{n\geq 0}$ such that $\deg{P_n}=n$ and  $\prodint{\un,P_n(x)P_m(x)}=K_{n}\delta_{n,m},$ where $\delta_{n,m}$ is the Kronecker  symbol and $K_{n}\ne 0$ (see \cite{Chihara,GMM21}). The sequence $(P_n)_{n\geq 0}$ is said to be the \textit{sequence of monic orthogonal polynomials} (SMOP, in short) with respect to the linear functional~$\un$. \\

Given   a quasi-definite linear functional  $\un$ and $(P_n)_{n\geq 0}$   its corresponding  SMOP,  there  exist two  sequences of complex numbers $( b_n)_{n\geq 0}$  and   $(a_n)_{n\geq 1}$, with $ a_n\ne 0$, such that	
\begin{equation}\label{ttrr}
	\begin{aligned}
		x\,P_{n}(x)&=P_{n+1}(x)+b_n\,P_{n}(x)+ a_{n}\,P_{n-1}(x),\quad n\geq 0,\\
		P_{-1}(x)&=0, \ \ \ \ \ P_{0}(x)=1,
	\end{aligned}
\end{equation}
where			
\begin{equation}\label{an}
	b_n=\frac{\langle\mathbf{u},x\,P^2_n(x)\rangle}{\langle\mathbf{u},P_n^2(x)\rangle}, \quad n\geq 0,\quad\quad a_n=\frac{\langle\mathbf{u},P_{n}^2(x)\rangle}{\langle\mathbf{u},P_{n-1}^2(x)\rangle},\quad n\geq 1 .
\end{equation}
If $\un$ is positive-definite, then $b_n$ and $a_n$ are real numbers with $a_n>0$.
Conversely, from Favard’s Theorem (see \cite{Chihara,GMM21}) if $(P_n)_{n\geq 0}$ is a sequence of monic polynomials  generated by a three term recurrence relation  as in \eqref{ttrr} with $ a_n\ne 0,$ $n\geq 1$, then there exists a unique linear functional $\un$ such that $(P_n)_{n\geq 0}$ is its SMOP. Moreover, in this case  the sequence of polynomials defined by  $$Q_{n}(x)=\dfrac{P_n(x)}{\prodint{\un,P_n^2}^{1/2}}$$ is said to be \textit{the sequence of orthonormal polynomials associated with $\un$}.

Another way to write the recurrence relation \eqref{ttrr} is in a matrix form. Indeed, if ${\bf P}=(P_0,P_1, \cdots)^\top,$  where $A^\top$ denotes the transposed of the matrix $A,$ then   \begin{equation}\label{matre}x{\bf P}=\J_{\operatorname{mon}}{\bf P},\end{equation} where $\J_{\operatorname{mon}}$ is the  semi-infinite matrix
\begin{equation*}
	\J_{\operatorname{mon}}=\begin{pmatrix}
		b_0 &1  &   &\\
		a_1 &b_1& 1&\\
		& a_2&b_2&\ddots\\
		& & \ddots&\ddots
	\end{pmatrix},
\end{equation*}
which is known in the literature as a \textit{monic Jacobi matrix} (see \cite{Chihara,GMM21}).
\begin{defi}\label{kernel}
	Let $(P_{n})_{n\ge 0}$ be   the   SMOP with respect to the quasi-definite moment functional~$\un$. We define the  $n$-th Christoffel-Darboux (C-D) kernel polynomial  as	$$K_n(x,y)=\sum_{k=0}^n\frac{P_k(x)P_k(y)}{\prodint{\un,P_k^2(x)}}.$$
\end{defi}

The classical \textit{(monic) Jacobi polynomials} are well-known families of polynomials depending on two parameters (see \cite{Chihara, OLB10}. Let $(P_n^{(\alpha,\beta)})_{n\geq0}$ be the monic Jacobi polynomials of parameters $(\alpha,\beta)$. For $\alpha,\, \beta>-1$, these polynomials are orthogonal with respect to the positive-definite linear functional $\un_{\alpha,\beta}$ defined by
$$
\prodint{\un_{\alpha,\,\beta},\,p(x)}=\int_{-1}^1 p(x)(1-x)^\alpha(1+x)^\beta dx,\quad p(x)\in\mathbb{P}.
$$
The explicit expression of these  polynomials is
\begin{equation*}
	P_n^{(\alpha,\beta)}(x)=\dfrac{1}{S_n(\alpha,\beta)}\sum_{k=0}^{n}\dbinom{n+\alpha}{n-k}\dbinom{n+\beta}{k}(x-1)^k(x+1)^{n-k},\ n\geq 0,\end{equation*}
where
$$ S_n(\alpha,\beta)=\dbinom{2n+\alpha+\beta}{n}.
$$ Here
$\binom{r}{k}=\frac{\Gamma(r+1)}{\Gamma(k+1)\Gamma(r-k+1)}$ is the usual binomial number and $$\Gamma(z)=\int_0^\infty t^{z-1}e^{-t}dt, \quad \operatorname{Re}(z)>0, $$ is the Gamma function. The explicit orthogonality relation reads
\begin{equation*}
	\prodint{\un_{\alpha,\,\beta},\,P_n^{(\alpha,\beta)}(x)P_m^{(\alpha,\beta)}(x)}={{2}^{2\,n+\alpha+\beta+1}\dfrac {\Gamma \left( n+\alpha+1 \right) \Gamma \left( n+\beta+
			1 \right) \Gamma \left( n+\alpha+\beta+1 \right)n! }{
			\left( 2\,n+\alpha+\beta+1 \right)  \left( \Gamma \left( 2\,n+\alpha+\beta+1 \right)
			\right) ^{2}}}\,
	\delta_{n,m}.
\end{equation*}
The Jacobi polynomials satisfy the three term recurrence relation
\begin{align*}
	&xP_n^{(\alpha,\beta)}(x)=P_{n+1}^{(\alpha,\beta)}(x)+b_n\,P_n^{(\alpha,\beta)}(x)+a_n\,P_{n-1}^{(\alpha,\beta)}(x), \quad n\geq 0,\\
	&P_{0}^{(\alpha,\beta)}(x)=1, \ \ P_{-1}^{(\alpha,\beta)}(x)=0,
\end{align*}
where
\begin{equation*}
	\begin{aligned}
		b_n&={\frac {{\beta}^{2}-{\alpha}^{2}}{ \left( 2\,n+2+\alpha+\beta \right)  \left( 2\,n+\alpha+\beta
				\right) }},\\
		\\
		a_n&=\,{\frac { 4\left( n+\beta \right)  \left( n+\alpha+\beta \right)  \left( n+\alpha
				\right) n}{ \left( 2\,n-1+\alpha+\beta \right)  \left( 2\,n+\alpha+\beta \right) ^{2}
				\left( 2\,n+\alpha+\beta+1 \right) }},
	\end{aligned}
\end{equation*}
except that when $\alpha=-\beta$, $b_0=\beta$ and $b_n=0$, $n\geq 1$.

We also have
\begin{equation*}
	P_n^{(\alpha,\beta)}(-x)=(-1)^nP_n^{(\beta,\alpha)}(x),\quad
	P_n^{(\alpha,\beta)}(1)=\dfrac{2^n}{S_n(\alpha,\beta)}\dbinom{n+\alpha}{n}.	
\end{equation*}
The representation of Jacobi polynomials in terms of a \textit{generalized Hypergeometric function} is
\begin{equation}\label{hiperrepre}
	P_n^{(\alpha,\beta)}(x)=\dfrac{2^n (\alpha+1)_n}{(\alpha+\beta+n+1)_n}\times\pFq{2}{1}\left({-n,\alpha+\beta+n+1};{\alpha+1};{\dfrac{1-x}{2}}\right),
\end{equation}
where $(a)_k$ is the \textit{Pochhammer symbol} defined by $(a)_0=1,$
\begin{equation}\label{poh}
	(a)_k=a(a+1)\cdots (a+k-1)=\dfrac{\Gamma(a+k)}{\Gamma(a)},	
\end{equation}
for $k=1,2,\ldots,$  and  \cite[pag. 352]{MMR94}
$$\pFq{p}{q}\left(a_1,\ldots,a_p;b_1,\ldots,b_q,x\right)=\sum_{k=0}^\infty\dfrac{(a_1)_k\cdots(a_p)_k}{(b_1)_k\cdots(b_q)_k}\dfrac{x^k}{k!}. $$
\subsection{Gauss quadrature formula} 
Given a continuous function $f(x)$ and a finite or infinite interval $I$, a problem of great interest is to approximate the integral
$$\int_I{f(x)}W(x)\,dx,$$
(where $W(x)$ is a weight function) by means of numerical integration  (quadrature rule)			
$$\int_I{f(x)}W(x)\,dx\approx f(x_1)w_1+\cdots +f(x_N)w_N.$$
Here,  the integral is approximated by a sum of  $N$ terms that involves the weights $w_k$ and the   nodes $x_k$ which must be properly chosen.
\begin{defi}
	We say that a quadrature rule has degree of exactness  $k$ if for every polynomial $q(x)\in\Px_k $,
	\begin{equation}\label{quad-rule}
		\int_I{q(x)}W(x)dx= q(x_1)w_1+\cdots q(x_n)w_N,
	\end{equation}				
	and there exists a polynomial of degree $k+1$ such that \eqref{quad-rule} does not hold.
\end{defi}
A simple way to construct quadrature rules is based on Lagrange interpolation formulas. Recall that given a set of $N$ points, $\mathcal{N}=\{(x_1,y_1),(x_2,y_2),\cdots, (x_N,y_N)\}$, where $x_i\neq x_j,$ if $i\ne j,$ the \textit{Lagrange interpolating polynomial} of degree at most $N-1$ corresponding to $\mathcal{N}$ is defined as
$$
L_N(x)=:\sum_{k=1}^N y_k\,\ell_k(x), \ \ \ \text{where}\ \ \ \  \ell_k(x)=:\prod_{\substack{ j=1\\ j\neq k}}^N\frac{(x-x_j)}{(x_k-x_j)}.
$$
Hence, when $N=1$, we take $L_1(x) \equiv y_1$. The numbers $x_j$, $1\le j \le N,$ are called \textit{nodes} and the polynomials  $(\ell_k(x))_{k=1}^N$ are known in the literature as the \textit{Lagrange polynomial basis}. The polynomial $L_N(x)$ has the property $L_N(x_j)=y_j$, $1\le j \le N$. Indeed, since $\ell_k(x_j)=\delta_{k,j}$, we have
$$
L_N(x_j)= \sum_{k=1}^N y_k\,\delta_{k,j}= y_j, \quad 1\le j \le N.
$$
The Lagrange polynomial $L_N(x)$ provides a unique solution to the problem of finding a polynomial of degree at most $N-1$ whose graph contains the points  $(x_j,y_j),$ $1\le j \le N$.

Given a continuous function  $f(x)$ and $N$ nodes $(x_j,f(x_j)),$ $1\le j \le N$, we can write
\begin{equation*}
	f(x)=L_N(x)+R_N(x),
\end{equation*}
where $R_N(x)$ is the interpolation error. Observe that if  $f(x)$ is a polynomial of degree at most $N-1,$  then by the uniqueness of Lagrange polynomial, $R_N(x)=0$.
With this in mind, we have the following quadrature rule 		
$$\int_I f(x)W(x)\,dx\approx\sum_{k=1}^N f(x_k)\int_I\ell_k(x)W(x)\,dx.$$
It is well known that in order to a quadrature formula has the maximum possible degree of exactness, the nodes $x_k$ must be chosen as the zeros of the orthogonal polynomials of degree $N$ with respect to the linear functional $\un$ associated with the weight function $W(x),$ i. e., $\prodint{\un, p}= \int_{I} p(x) W(x) dx.$ Let remind that such zeros are real, simple and are located in the interior of the convex hull of the support of the measure associated with the linear functional (see \cite{Chihara}). In other words:
\begin{teo}[Gauss Quadrature Rule]
	Let $\un$ be a positive-definite linear  functional and let $(P_n)_{n\ge 0}$ be its corresponding SMOP. There exist positive real numbers $A_{N,1},\ldots, A_{N,N}$ such that for every polynomial   $q(x)$ of degree at most  $2N-1$, we have				
	\begin{equation*}
		\prodint{\un,q(x)}=\sum_{k=1}^NA_{N,k}\,q(x_{N,k}),
	\end{equation*}
	where $(x_{N,k})_{k=1}^N$ are the zeros of  de $P_N(x)$.
\end{teo} 	
The constants $A_{N,k}$ involved in the Gauss quadrature rule are known in the literature as \textit{Christoffel numbers}.
\begin{coro}
	The Christoffel numbers  $A_{N,k},$ $1\le k \le N$, have the following representation,
	\begin{equation}\label{Chistoffelnumbers}
		A_{N,k}=\frac{1}{K_{N-1}(x_{N,k},x_{N,k})},	
	\end{equation}
	where $K_{N-1}(x,y)$ is the $(N-1)$-th C-D kernel polynomial as in Definition \ref{kernel}.
\end{coro}
Notice that from \eqref{matre} we get the relation
$$x_{N,k}\begin{pmatrix}P_0(x_{N,k})\\\vdots\\P_{N-1}(x_{N,k})\end{pmatrix}=(\J_{\operatorname{mon}})_{N\times N}\begin{pmatrix}P_0(x_{N,k})\\\vdots\\P_{N-1}(x_{N,k})\end{pmatrix},$$
where $(\J_{\operatorname{mon}})_{N\times N}$ is the $N\times N$ leading principal submatrix of $\J_{\operatorname{mon}}$. This means that $x_{N,k}$ is an eigenvalue of $(\J_{\operatorname{mon}})_{N\times N}$ with associated eigenvector  $(P_0(x_{N,k}),\cdots,P_{N-1}(x_{N,k}))^T$.  The matrix $(\J_{\operatorname{mon}})_{N\times N}$ can be 	symmetrized by a diagonal similarity transformation $\D$ to obtain
$$\J_{N\times N}=\mathbf{D}(\J_{\operatorname{mon}})_{N\times N}\mathbf{D}^{-1}=\begin{pmatrix}
	b_0&\sqrt{a_1}&&&\\[5pt]
	\sqrt{a_1}&b_1&\sqrt{a_2}&&\\
	&\ddots&\ddots&\ddots&\\
	&&\sqrt{a_{N-2}}&b_{N-2}&\sqrt{a_{N-1}}\\[5pt]
	&&&\sqrt{a_{N-1}}&b_{N-1}\\
\end{pmatrix}.$$
Notice that $\J_{N\times N}$ is the $N\times N$ leading principal submatrix of the Jacobi matrix associated with the orthonormal polynomials.  Since the eigenvalues are preserved
by a similarity transformation, $x_{N,k}$ is also an eigenvalue of the symmetric tridiagonal matrix $\J_{N\times N}$.
Moreover, taking into account \eqref{Chistoffelnumbers}, if $v_k$ is the eigenvector corresponding to $x_{N,k}$ normalized as $\|v_{k}\|_2 = 1$, then (see \cite[Theorem 3.1]{M15})
$$A_{N,k}=\un_0\,v^2_{k,1},$$ where $v_{k,1}$ is the first component of $v_k$.  It is well known   that the seeking of eigenvalues and eigenvectors of a symmetric tridiagonal matrix is a relatively efficient and well-conditioned procedure.

Notice that to use the above method we  need to know the explicit values of $(b_k)_{k=0}^{N-1}$ and $(a_k)_{k=1}^{N-1}$,
but it is not always easy from \eqref{an}. Moreover, the usual  approach is to represent each $P_k(x)$ explicitly as a polynomial in $x$ and
to compute the inner products by multiplying out term by term. This will be feasible if we
know the first $2N$ moments
\begin{equation}\label{momentos}
	\un_k=\int_{I} x^k W(x)\,dx, \quad k=0,\ldots,2N-1.
\end{equation}
However, the solution of the resulting set of algebraic equations for the coefficients $b_k$ and $a_k$
in terms of the moments $\un_k$ is in general  ill-conditioned. Sack and Donovan \cite{SD71} (see also \cite{G70}) showed  that the numerical stability is improved if
instead of using the classical moments $\un_k$  we use the modified moments
\begin{equation}\label{modfmomentos}
	m_k =\int_{I} \pi_k(x)W(x)\,dx,
\end{equation}
where $(\pi_k(x))_{k=0}^{2N-1} $ is a set of $2N$ monic polynomials satisfying a three-term recurrence relation
\begin{equation}\label{cofmod}
	\begin{aligned}
		x\pi_{k}(x) &= \pi_{k+1}(x)+{\mathfrak b}_k\pi_k(x)+{\mathfrak a_k}\pi_{k-1}(x), \ \ k\geq 0,\\
		\pi_0(x) &= 1,\ \  \pi_{-1}(x) = 0,	
	\end{aligned}	
\end{equation}
with known coefficients $\mathfrak{b}_k\in\mathbb{R}$ and $\mathfrak{a}_k\geq 0$. Notice that when $\mathfrak{b}_k = \mathfrak{a}_k = 0$, for all $k \geq 0$, the modified  moments \eqref{modfmomentos} become the standard moments \eqref{momentos}. Thus, if we get the modified moments \eqref{modfmomentos}, then we can find the coefficient as follows:\\

Assume that we know the modified moments
$$m_k =\int_{I} \pi_k(x)W(x)dx,\quad j=0,\ldots, 2N-1,$$
where the $\pi_k(x)$’s satisfy  the recurrence relation \eqref{cofmod}
and the coefficients $\mathfrak{b}_k$, $\mathfrak{a}_k$ are known explicitly. In Sack and Donovan \cite{SD71} an efficient algorithm to determine $b_k$ and $a_k$ via a set of intermediate quantities is given. The  explanation of the following algorithm can be found in \cite{Ga78,W74}.
$$\sigma_{k,\ell}=\prodint{\un,p_k(x)\pi_\ell(x)},\quad k,\ \ell\geq -1,$$
where $$\sigma_{-1,\ell}=0,\quad \ell=1,2,\ldots,2N-2,$$
$$\sigma_{0,\ell}=m_\ell,\quad \ell=0,1,\ldots,2N-1,$$
$$b_0=\mathfrak{b_0}+\dfrac{m_1}{m_0},\qquad a_0=0,$$
and for $k=1,2,\ldots,N-1,$
\begin{equation}\label{chevalgo}
	\begin{aligned}
		\sigma_{k,\ell}&=\sigma_{k-1,\ell+1}-(b_{k-1}-\mathfrak{b}_{\ell})\,\sigma_{k-1,\ell}-a_{k-1}\,\sigma_{k-2,\ell}+\mathfrak{a}_{\ell}\,\sigma_{k-1,\ell-1},\quad \ell=k,\ldots, 2N-k-1,\\
		b_k&=\mathfrak{b}_k-\dfrac{\sigma_{k-1,k}}{\sigma_{k-1,k-1}}+\dfrac{\sigma_{k,k+1}}{\sigma_{k,k}},	\\
		a_k&=\dfrac{\sigma_{k,k}}{\sigma_{k-1,k-1}}.
	\end{aligned}
\end{equation}
The above algorithm is said to be \textit{modified Chebyshev algorithm.}
\section{Quadrature Formula for the truncated Laguerre polynomials}\label{S3}
Let $(L^\alpha_{n} (x;z))_{n\geq0}$ be the sequence of monic orthogonal polynomials with respect to the linear functional
\begin{equation}\label{laguerrefunc}
	\langle\elln,p(x)\rangle=\int_0^{z} p(x)\,x^{\alpha} e^{-x}dx, \quad p(x)\in\mathbb{P},
\end{equation}
with $\alpha>-1$ and $z>0$.
This functional is known in the literature as \textit{the truncated Gamma functional} and the polynomials $ (L^\alpha_{n} (x;z))_{n\geq 0}$ are called \textit{truncated Laguerre orthogonal polynomials} (see \cite{DGM23}).
Notice that the moments of the linear functional $\elln$ are given by
$$\elln_m=\int_0^z x^{m+\alpha}e^{-x}dx=\widehat\gamma(m+\alpha+1,z),$$
where $\widehat\gamma(a,z)$ is the \textit{incomplete gamma function} defined by \cite[8.2.1]{OLB10}
$$\widehat\gamma(a,z)=\int_0^z x^{a-1}e^{-x}dx. $$ Its  series representation, (see \cite[8.5.1]{OLB10}), is
$$\widehat\gamma(a,z)=\dfrac{z^ae^{-z}}{a}\sum_{k=0}^\infty\dfrac{z^k}{k(a+1)_k}.$$
A change of variables in \eqref{laguerrefunc} yields
$$\prodint{\elln,p}=z^\alpha\int_0^1p(zx) x^\alpha e^{-zx}dx.$$

According to \cite[Eq. 7.10--7.11]{DGM23} we get that the recurrence coefficients associated with  the monic orthogonal polynomials \eqref{ttrr}  with respect to the weigh $w(x,z)=x^\alpha e^{-zx}$ supported on $[0,1]$ satisfy	
\begin{equation*}
	\begin{aligned}
		a_{n}\left(  z\right)  =&\dfrac{1}{16}+\dfrac{\left(  1-2\alpha^{2}\right)}{64}n^{-2}+\alpha \dfrac{\alpha(z+2\alpha)-1}{64}n^{-3}+\mathcal{O}\left(n^{-4}\right)  , 
		\\
		b_{n}\left( z\right)  =&\dfrac{1}{2}-\dfrac{\alpha^{2}}{8}n^{-2}+\dfrac{z-2\alpha^{2}(z+2\alpha+2)}{32}n^{-3}+\mathcal{O}\left(  n^{-4}\right)  ,
	\end{aligned}
\end{equation*}
as $n\to \infty$, or, equivalently,
\begin{equation}\label{limit}
	\lim_{n\to \infty}a_n(z)=\dfrac{1}{16},\qquad \lim_{n\to \infty}b_n(z)=\dfrac{1}{2}.	
\end{equation}
Notice also that  $w(x,z)\to x^{\alpha}$ when $z\to 0$. Therefore, a good choice of sequence of monic orthogonal polynomials to use in the modified Chebyshev algorithm is the \textit{shifted Jacobi polynomials} with parameters $(0,\alpha)$, which are orthogonal with respect to the weight function $w(x)=x^\alpha$ supported  on the interval $[0,1]$.
\begin{pro}
	The shifted	monic Jacobi polynomials $(\widetilde P^{(0,\alpha)}_n)_{n\geq 0}$ can be written as
	\begin{equation}\label{repre1}
		\widetilde{P}^{(0,\alpha)}_n(x)=\dfrac{1}{2^n}P^{(0,\alpha)}_n(2x-1)=\dfrac{(-1)^{n}}{2^n}P^{(\alpha,0)}_n(1-2x).	
	\end{equation}
	where $P^{(0,\alpha)}_n(x)$ is the Jacobi polynomial of parameters $(0,\alpha)$ and degree $n$.
\end{pro}

\begin{pro}
	Let $(\wtP^{(0,\alpha)}_{n})_{n\geq 0}$ be the sequence of  monic shifted Jacobi polynomials orthogonal with respect to  the weight function
	$\omega(x)=x^\alpha$, $\alpha>-1$, on $[0,1]$. The modified moments
	$$\begin{aligned}
		m_n(\alpha,z)=\int_{0}^1 \wtP_n^{(0,\alpha)}(x)x^{\alpha}e^{-zx}dx,	
	\end{aligned}$$
	where $n = 0, 1,\ldots,$ can be expressed via a $(1,1)$-Generalized Hypergeometric function as
	$$\begin{aligned}
		m_0(\alpha,z)&=\begin{cases}
			\dfrac{1}{\alpha+1},& z=0,\\[10pt]
			\dfrac{\widehat\gamma(\alpha+1,z)}{z^{\alpha+1}},& z>0,		
		\end{cases}\\
		m_n(\alpha,z)&=(-1)^n \dfrac{n!}{(\alpha+n+1)_n^2(\alpha+2n+1)}z^ne^{-z}\times\pFq{1}{1}\left(n+1;\alpha+2n+2;z\right),\quad n\geq 1.
	\end{aligned}$$
\end{pro} 	
\begin{proof}
	From \cite[eq. 7.523]{GR80},
	\begin{multline*}
		\int_{0}^1 x^{v-1}(1-x)^{\rho-1}e^{-zx} \times\pFq{2}{1}\left(r,q;v;x\right)dx=\\ \dfrac{\Gamma(v)\Gamma(\rho)\Gamma(v+\rho-r-q)}{\Gamma(v+\rho-r)\Gamma(v+\rho-q)}e^{-z}\times\pFq{2}{2}\left(\rho,v+\rho-r-q;v+\rho-r,v+\rho-q;z\right)
	\end{multline*}
	holds for
	$$\operatorname{Re}(v)>0,\quad \operatorname{Re}(\rho)>0,\quad\text{and}\quad  \operatorname{Re}(v+\rho-r-q)>0.$$
	In particular, for $\varepsilon>0$ and $n\geq 0$
	\begin{equation}\label{inth1}
		\begin{aligned}
			&\int_{0}^1 x^{(\alpha+1)-1}e^{-xz} \,\pFq{2}{1}\left(-n-\varepsilon,\alpha+n+\varepsilon+1;\alpha+1;x\right)dx\\ &=\dfrac{\Gamma(\alpha+1)}{\Gamma(\alpha+2+n+\varepsilon)\Gamma(1-n-\varepsilon)}e^{-z} \times\pFq{2}{2}\left(1,1;\alpha+2+n+\varepsilon,1-n-\varepsilon;z\right)\\
			&=\dfrac{\Gamma(\alpha+1)}{\Gamma(\alpha+2+n+\varepsilon)\Gamma(1-n-\varepsilon)}e^{-z}\sum_{k=0}^\infty\dfrac{(1)_k^2}{(\alpha+n+\varepsilon+2)_k(1-n-\varepsilon)_k}\dfrac{z^k}{k!}\\
			&=	\Gamma(\alpha+1)e^{-z}\sum_{k=0}^\infty\dfrac{(1)_k^2}{\Gamma(\alpha+n+\varepsilon+2+k)\Gamma(1-n-\varepsilon+k)}\dfrac{z^k}{k!}.
		\end{aligned}
	\end{equation}
	On the other hand,  from \eqref{hiperrepre} we get
	\begin{align*}
		\dfrac{(-1)^n}{2^n}P^{(\beta,\alpha)}_n(-x)=(-1)^n\dfrac{ (\beta+1)_n}{(\alpha+\beta+n+1)_n}\times\pFq{2}{1}\left({-n,\alpha+\beta+n+1};{\beta+1};{\dfrac{1+x}{2}}\right),
	\end{align*}
	and by  \eqref{repre1}
	\begin{align*}
		\wtP^{(0,\alpha)}_{n}(x)&=(-1)^nP^{(\alpha,0)}_n(-(2x-1))\\
		&=(-1)^n\dfrac{ (\alpha+1)_n}{(\alpha+n+1)_n}\times\pFq{2}{1}\left({-n,\alpha+n+1};{\alpha+1};x \right).
	\end{align*}
	Thus, from \eqref{inth1} we obtain
	\begin{align*}
		&(-1)^n\dfrac{(\alpha+n+1)_n}{ (\alpha+1)_n}\int_{0}^1 \wtP_n^{(0,\alpha)}(x)x^{\alpha}e^{-zx}dx\\ &=\int_{0}^1 x^{(\alpha+1)-1}e^{-zx} \,\lim_{\varepsilon\to 0}\,\pFq{2}{1}\left(-n-\varepsilon,\alpha+n+1+\varepsilon;\alpha+1;x\right)dx\\
		\\ &=\lim_{\varepsilon\to 0}\int_{0}^1 x^{(\alpha+1)-1}e^{-zx} \,\pFq{2}{1}\left(-n-\varepsilon,\alpha+n+1+\varepsilon;\alpha+1;x\right)dx\\
		&=\Gamma(\alpha+1)e^{-z}\lim_{\varepsilon\to 0}	\sum_{k=0}^\infty\dfrac{(1)_k^2}{\Gamma(\alpha+n+\varepsilon+2+k)\Gamma(1-n-\varepsilon+k)}\dfrac{z^k}{k!}.
	\end{align*}
	Notice that $\dfrac{1}{\Gamma(1-n-\epsilon+k)}\to 0$ when $\varepsilon\to 0$ for $k=0,\ldots, n-1.$ Moreover,  using  Stirling's Formula (see \cite[Pag. 141, Eq. 5.11.7]{OLB10}) for $k\gg 0$ we get
	\begin{multline*}
		\left|\dfrac{(1)_k^2}{\Gamma(\alpha+n+\varepsilon+2+k)\Gamma(1-n-\varepsilon+k)}\dfrac{z^k}{k!}\right|=	\left|\dfrac{k!}{(\alpha+n+\varepsilon+2)_k\Gamma(1-n-\varepsilon+k)}z^k\right|\approx\\[10pt]
		\left|\dfrac{\sqrt{ k}\left(\dfrac{k}{e}\right)^{k}}{\sqrt{2\pi(\alpha+n+\varepsilon+1+k)(k-n-\varepsilon)}\left(\dfrac{\alpha+n+\varepsilon+1+k}{e}\right)^{\alpha+n+\varepsilon+1+k}\left(\dfrac{k-n-\varepsilon}{e}\right)^{k-n-\varepsilon}}\,z^k\right|.
	\end{multline*}
	If we denote by $s_k$ the last expression, then we get
	$$\lim_{k\to \infty}(s_{k})^{1/k}=0.$$  Consequently, from Tannery's theorem (see \cite[page 448]{MK05}) we can interchange the limit and the summatory. Thus, using   \eqref{poh} we get
	\begin{align*}
		(-1)^n\dfrac{(\alpha+n+1)_n}{(\alpha+1)_n}\int_{0}^1 \wtP_n^{(0,\alpha)}(x)x^{\alpha}e^{-zx}dx
		&=\Gamma(\alpha+1)e^{-z}\sum_{k=n}^\infty\dfrac{(1)^2_k}{\Gamma(\alpha+2+n+k)\Gamma(k+1-n)}\dfrac{z^k}{k!}\\
		&=\dfrac{\Gamma(\alpha+1)\,n!}{\Gamma(\alpha+2n+2)}z^ne^{-z}\sum_{t=0}^\infty\dfrac{(n+1)_t}{(\alpha+2n+2)_{t}}\dfrac{z^{t}}{t!}
		\\
		&=\dfrac{\Gamma(\alpha+1)\,n!}{\Gamma(\alpha+2n+2)}z^ne^{-z}\times\pFq{1}{1}\left(n+1;\alpha+2n+2;z\right).
	\end{align*}
\end{proof}
\section{Numerical experiments}\label{S4}
By default, \texttt{MATLAB} works with  16 digits of precision. For higher precision, we can use the function  \texttt{vpa(expr,digits)}  in \texttt{Symbolic Math Toolbox}, where \texttt{expr} is the numerical variable and \texttt{digits} is the number of   precision digits which can be increased from 1 until $2^{29}+1$. By default, the variable-precision arithmetic of  \texttt{vpa} uses $32$ significant decimal digits of precision. For our propose we take \texttt{digits=100} and we assume that the values obtained with this  precision are "exact".\\

In this section we will use the \textit{Modified Chebyshev algorithm and the quadrature algorithm}   which are collected in the package \texttt{OPQ} (Orthogonal Polynomials and Quadrature) in MATLAB provided by Walter Gautschi. All pieces of software into  the package \texttt{OPQ}  can be downloaded from \cite{G18}
\begin{center}
	\url{www.siam.org/books/se26} \end{center}
The explanation of these routines is illustrated in \cite{G05,G16,G83}. In particular,  the modified Chebyshev moment   routine \eqref{chevalgo}  for generating the  $N$   first coefficients  $b_n$ and  $a_n$ associated with  the polynomials $(L_n^{\alpha}(x;z))_{n= 0}^{N-1}$  we use:
\begin{center}
	{\texttt ab = chebyshev(2N, mom, abm),	}
\end{center}
where $\texttt{mom}=(m_0(\alpha,z),m_1(\alpha,z),\cdots, m_{2N-1}(\alpha,z))$ is the $1\times 2N$ array of modified moments and \texttt{abm} the
$(2N-1)\times 2$ array of known recurrence coefficients
$$\texttt{abm} =\begin{pmatrix}
	\mathfrak{b}_0 &\mathfrak{a}_0\\
	\mathfrak{b}_1 &\mathfrak{a}_1\\
	\vdots& \vdots\\
	\mathfrak{b}_{2N-2}& \mathfrak{a}_{2N-2}
\end{pmatrix}.$$
The corresponding symbolic routine is provided in the package \texttt{SOPQ}  \cite{G18}
\begin{center}
	{\texttt ab = schebyshev(digits, 2N, mom, abm).}
\end{center}	
Taking into account the above, we define \textit{the maximal relative error} by 	
\begin{equation}\label{error}
	\mathcal{E}_n^\alpha(z;\operatorname{digits}):=\max_{0\leq k\leq n-1}\left\{\left|\dfrac{a_k-\mathrm{\bf a}_k}{\mathrm{\bf a}_k}\right|,\ \left|\dfrac{b_k-\mathrm{\bf b}_k}{\mathrm{\bf b}_k}\right|\right\},\end{equation}
where  the "exact values" $\mathrm{\bf a}_k$ and $\mathrm{\bf b}_k$ are obtained using the symbolic routine	
\begin{center}
	\texttt{schebyshev(digits, 2N, mom, abm)}.	
\end{center}
with \texttt{digits=100}. For example, for $\alpha=1$, $z=1,$ $N=48$ and \texttt{digits=30}, we get the first 48 recurrence   coefficients $(b_k)_{k=0}^{47}$, $(a_k)_{k=0}^{47}$ (see Table \ref{tabla1} and Table \ref{tabla2}).

\begin{table}[ht]
	\begin{equation*}
		\renewcommand{\arraystretch}{1.1}\begin{array}{|c|c|c|c|c|c|}\hline
			k&b_k(1)&k&b_k(1)&k&b_k(1)\\ \hline\hline
			0&   0.607788808822667&  16&  0.500426808722825&32&  0.500113960288728\\
			1&   0.531655773460623&  17&  0.500380956440347&33&  0.500107376935035\\
			2&   0.513427891887575&  18&  0.500342117324710&34&  0.500101348056586\\
			3&   0.507534434505874&  19&  0.500308930359149&35&  0.500095813084832\\
			4&   0.504835172165686&  20&  0.500280349444746&36&  0.500090719501977\\
			5&   0.503368518898357&  21&  0.500255559818526&37&  0.500086021590060\\
			6&   0.502482032582844&  22&  0.500233919258297&38&  0.500081679401089\\
			7&   0.501905070635327&  23&  0.500214916014164&39&  0.500077657904668\\
			8&   0.501508458210362&  24&  0.500198138240255&40&  0.500073926279007\\
			9&   0.501224082197933&  25&  0.500183251468526&41&  0.500070457318386\\
			10&  0.501013226345908&  26&  0.500169981794040&42&  0.500067226935724\\
			11&  0.500852549532634&  27&  0.500158103174469&43&  0.500064213743185\\
			12&  0.500727295329595&  28&  0.500147427732113&44&  0.500061398697152\\
			13&  0.500627761645442&  29&  0.500137798273680&45&  0.500058764796504\\
			14&  0.500547355966703&  30&  0.500129082466506&46&  0.500056296825256\\
			15&  0.500481471186175&  31&  0.500121168264855&47&  0.500053981132248\\ \hline
		\end{array}
	\end{equation*}
	\caption{Sequence  $b_k(1)$, for $k=0,\ldots,47.$ }
	\label{tabla1}
\end{table}
\begin{table}[ht]
	\begin{equation*}
		\renewcommand{\arraystretch}{1.1}\begin{array}{|c|c|c|c|c|c|}\hline
			k&a_k(1)&k&a_k(1)&k&a_k(1)\\ \hline		\hline
			0&      0.2642411176571153 &  16&   0.06244605779295135 &32&  0.06248566001620657  \\
			1&     0.06174799916059207 &  17&   0.06245187208765875 &33&  0.06248649068343801  \\
			2&     0.06110639779446930 &  18&   0.06245679527042237 &34&  0.06248725122345314  \\
			3&     0.06159770170459388 &  19&   0.06246100042276222 &35&  0.06248794931205289  \\
			4&     0.06190038948667886 &  20&   0.06246462058754074 &36&  0.06248859160306462  \\
			5&     0.06207694279947663 &  21&   0.06246775938221266 &37&  0.06248918388734983  \\
			6&     0.06218661764444579 &  22&   0.06247049847158962 &38&  0.06248973122368508  \\
			7&     0.06225891009402938 &  23&   0.06247290291817096 &39&  0.06249023804706056  \\
			8&     0.06230892605966116 &  24&   0.06247502507158938 &40&  0.06249070825874153  \\
			9&     0.06234491040894508 &  25&   0.06247690743561846 &41&  0.06249114530151897  \\
			10&    0.06237163970905703 &  26&   0.06247858480860449 &42&  0.06249155222286876  \\
			11&    0.06239202635370342 &  27&   0.06248008590029081 &43&  0.06249193172818968  \\
			12&    0.06240792454776711 &  28&   0.06248143456640384 &44&  0.06249228622586166  \\
			13&    0.06242055888438060 &  29&   0.06248265076085398 &45&  0.06249261786552981  \\
			14&    0.06243076389952361 &  30&   0.06248375127700079 &46&  0.06249292857075320  \\
			15&    0.06243912391137471 &  31&   0.06248475032972568 &47&  0.06249322006694707  \\ \hline		
		\end{array}
	\end{equation*}
	\caption{Sequence  $a_k(1)$, for $k=0,\ldots,47.$}
	\label{tabla2}
\end{table}

Notice that from here we infer that the convergence of the previous sequences to the limits \eqref{limit} is satisfied. However,  it convergence appears to be "slow" since for  $k=47$ we only have $4$ significant digits.

In order to study the stability of our procedure, we will take $\alpha=1$ and we will  use \eqref{error} with two different arithmetic, the usual \texttt{digits=16} and the high precision  \texttt{digits=30} for different values of~$z$.
\begin{table}[ht]
	\renewcommand{\arraystretch}{1.2}\begin{tabular}{|c| c|c|c|}\hline
		\multicolumn{4}{|c|}{$\alpha=1$}\\ \hline
		\texttt{digits} &$z=5$ & $z=10$& $z=15$  \\    \hline\hline
		16 &$7.8056\times10^{-16}$ & $  8.8163\times10^{-15}$ & $6.0238\times10^{-13}$\\ \hline
		30 &$7.8238\times10^{-38}$ & $8.3718\times10^{-37}$   & $4.9234\times10^{-35}$ \\ \hline
	\end{tabular}	
	\caption{Maximal relative error $\mathcal{E}_{50}^1(z;\operatorname{digits})$ for $z=5,\,10,\,15$.}
\end{table}
\begin{table}[ht]
	\renewcommand{\arraystretch}{1.2}\begin{tabular}{|c|c|c|c|}\hline
		\multicolumn{4}{|c|}{$\alpha=1$}\\ \hline
		\texttt{digit}& $z=20$ &$z= 25$ & $z=30$\\ \hline\hline
		16    &$ 2.6790\times10^{-11}$ &$4.8261\times10^{-9}$ & $4.2894\times 10^{-7}$\\ \hline
		30   & $3.6960\times10^{-33}$ &$5.9650\times10^{-31}$ & $9.9114\times10^{-30}$\\ \hline
	\end{tabular}	
	\caption{Maximal relative error $\mathcal{E}_{50}^1(z;\operatorname{digits})$ for $z=20,\,25,\,30$. }
\end{table}
The loss of digits practically does not exist for  any  value of $z\in[0,30]$ and \texttt{digits=30}. For \texttt{digits=16} the loss of digits practically does not exist for  any  value of $z\in[0,10)$. When $z=10,$ $z=15,$ $z=20,$ $z=25,$ and $z=30,$  we have a loss of at most one, three, five, seven or nine  digits, respectively, regardless of which arithmetic we use. Thus, for $z < 10$, the method is well-conditioned and its condition number
is near $1$. On the other hand, taking into account \cite{M15}, if the condition number of a mapping is $10^t$, then approximately $t$ decimal digits are lost. Therefore, we can conclude  that the condition
number is approximately $10$, $10^3$, $10^5$ ,$10^7$  and $10^9$ for $z=10,$ $z=15,$ $z=20,$ $z=25,$ and $z=30,$ respectively. Thus, from the above, we conclude that if we want a  accuracy of $d$ decimal digits in the recurrence coefficients, then we must work with \texttt{digit}$=t+d,$ \cite{M15}.\\

Suppose now, that  we want to determinate the coefficients  $b_k$ and $a_k$ for $k=0,\ldots,n-1$, where by hypothesis $\mathfrak{a}_0=m_0(\alpha,z)$. For simplicity assume that   $\alpha=1$ and $z\in [0,T]$ with $T>0$. Taking into account  \cite{S15}, we first select a system of points $Q = \{z_k\}$ in the interval $[0,T]$.  In these points we determine the corresponding values of the coefficients $(b_k(z))_{k=0}^{n-1}$, $(a_k(z))_{k=1}^{n-1}$ and then construct the corresponding interpolating functions for each one of these coefficients, using the function \texttt{interp1(Q,f\_Q,q,"spline")}, where \texttt{f\_Q} is the vector of values of the function~$f$ corresponding to the partition   $Q$ and $q$ is the interval where the interpolator function is defined. In this case, the interpolation is made using spline's technique.

\begin{figure}[ht]
	\includegraphics[width=13cm]{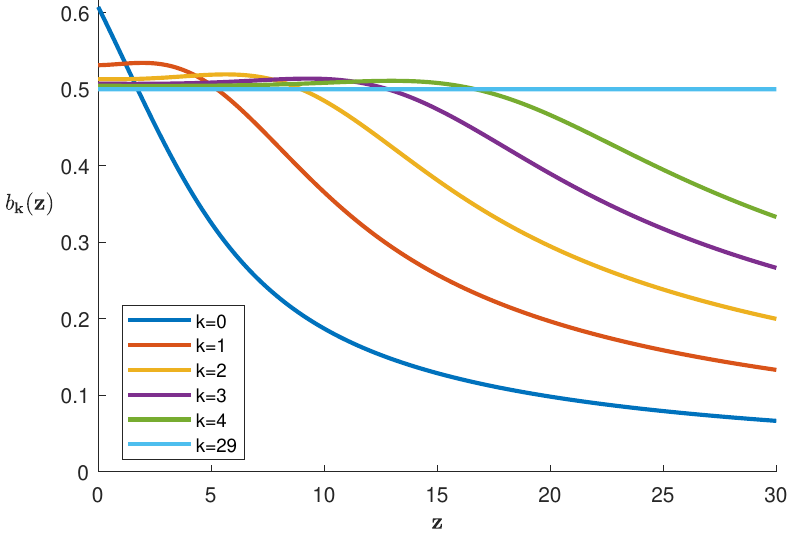}
	\caption{The coefficients $b_k(z)$ for $z\in[0,30]$ in the case $\alpha=1$.}
	\label{gra-b-k}
\end{figure}

\begin{figure}[ht]
	\includegraphics[width=13cm]{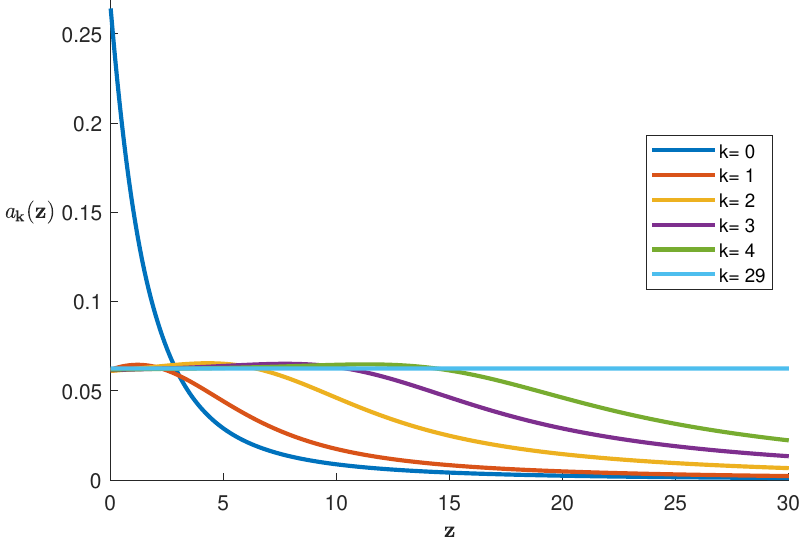}
	\caption{The coefficients $a_k(z)$ for $z\in[0,30]$ in the case $\alpha=1$.}
	\label{gra-a-k}
\end{figure}

We illustrate this procedure for $\alpha=1,$ $T=30$. By simplicity we take an equidistant partition of points
$$Q=\{z_{k+1}=z_{k}+3/4, \quad z_0=0, \quad k=0,\ldots, 29\},$$
or,  equivalently in \texttt{MATLAB} \texttt{linspace(0,30,90)}. In Figure \ref{gra-b-k} and Figure \ref{gra-a-k} we present the coefficient $b_k(z),$ $a_k(z)$ for $\alpha=1.$

Finally, we use the function \texttt{gauss(N,ab)} from  package \texttt{OPQ} to obtain the zeros of the polynomial $L^{\alpha}_{50}(x; z)$ (in increasing order) as well as the corresponding Christoffel numbers~\eqref{Chistoffelnumbers}. If we take  $N=50$, $\alpha=1$ and $z=1$ we get Table~\ref{gauscuadrature}.

\begin{table}[ht]
	\begin{equation*}	
		\renewcommand{\arraystretch}{1.1}\begin{array}{|c|c|} \hline
			\lambda_{50,k}&A_{50,k}\\	\hline\hline
			1.0723\times 10^{-3} &  1.8699\times 10^{-6}\\
			3.5934\times 10^{-3} &  1.0454\times 10^{-5}\\
			7.5515\times 10^{-3} &  2.8145\times 10^{-5}\\
			1.2941\times 10^{-2} &  5.3562\times 10^{-5}\\
			1.9753\times 10^{-2} &  8.2141\times 10^{-5}\\
			2.7979\times 10^{-2} &  1.0792\times 10^{-4}\\
			3.7607\times 10^{-2} &  1.2564\times 10^{-4}\\
			4.8621\times 10^{-2} &  1.3232\times 10^{-4}\\
			6.1005\times 10^{-2} &  1.2782\times 10^{-4}\\
			7.4740\times 10^{-2} &  1.1435\times 10^{-4}\\
			8.9804\times 10^{-2} &  9.5435\times 10^{-5}\\
			1.0617\times 10^{-1} &  7.4720\times 10^{-5}\\
			1.2382\times 10^{-1} &  5.5130\times 10^{-5}\\
			1.4270\times 10^{-1} &  3.8478\times 10^{-5}\\
			1.6281\times 10^{-1} &  2.5488\times 10^{-5}\\
			1.8408\times 10^{-1} &  1.6071\times 10^{-5}\\
			2.0648\times 10^{-1} &  9.6726\times 10^{-6}\\\hline
		\end{array}
		\renewcommand{\arraystretch}{1.1}\begin{array}{|c|c|} \hline
			\lambda_{50,k}&A_{50,k}\\ \hline\hline
			2.2997\times 10^{-1} &  5.5708\times 10^{-6}\\
			2.5450\times 10^{-1} &  3.0779\times 10^{-6}\\
			2.8001\times 10^{-1} &  1.6354\times 10^{-6}\\
			3.0644\times 10^{-1} &  8.3766\times 10^{-7}\\
			3.3373\times 10^{-1} &  4.1466\times 10^{-7}\\
			3.6181\times 10^{-1} &  1.9888\times 10^{-7}\\
			3.9059\times 10^{-1} &  9.2662\times 10^{-8}\\
			4.2001\times 10^{-1} &  4.2055\times 10^{-8} \\
			4.4998\times 10^{-1} &  1.8645\times 10^{-8}\\
			4.8039\times 10^{-1} &  8.0995\times 10^{-9}\\
			5.1114\times 10^{-1} &  3.4581\times 10^{-9}\\
			5.4214\times 10^{-1} &  1.4559\times 10^{-9}\\
			5.7326\times 10^{-1} &  6.0646\times 10^{-10}\\
			6.0437\times 10^{-1} &  2.5087\times 10^{-10}\\
			6.3536\times 10^{-1} &  1.0343\times 10^{-10}\\
			6.6609\times 10^{-1} &  4.2669\times 10^{-11}\\
			6.9640\times 10^{-1} &  1.7684\times 10^{-11}\\\hline
		\end{array}
		\renewcommand{\arraystretch}{1.1}\begin{array}{|c|c|} \hline
			\lambda_{50,k}&A_{50,k}\\\hline\hline
			7.2616\times 10^{-1} &  7.3924\times 10^{-12}\\
			7.5520\times 10^{-1} &  3.1302\times 10^{-12}\\
			7.8337\times 10^{-1} &  1.3481\times 10^{-12}\\
			8.1051\times 10^{-1} &  5.9302\times 10^{-13}\\
			8.3645\times 10^{-1} &  2.6750\times 10^{-13}\\
			8.6103\times 10^{-1} &  1.2420\times 10^{-13}\\
			8.8407\times 10^{-1} &  5.9564\times 10^{-14}\\
			9.0543\times 10^{-1} &  2.9587\times 10^{-14}\\
			9.2495\times 10^{-1} &  1.5253\times 10^{-14}\\
			9.4247\times 10^{-1} &  8.1661\times 10^{-15}\\
			9.5787\times 10^{-1} &  4.5341\times 10^{-15}\\
			9.7102\times 10^{-1} &  2.5981\times 10^{-15}\\
			9.8182\times 10^{-1} &  1.5189\times 10^{-15}\\
			9.9016\times 10^{-1} &  8.8292\times 10^{-16}\\
			9.9598\times 10^{-1} &  4.7777\times 10^{-16}\\
			9.9924\times 10^{-1} &  1.8745\times 10^{-16}\\
			&\\
			\hline
		\end{array}
	\end{equation*}
	\caption{Nodes and weight coefficients. Gauss Quadrature formula}
	\label{gauscuadrature}
\end{table}
\section{Concluding remarks}
The  Gauss-Rys quadrature formula with respect to the weight function $w(x)=e^{-tx^2}$ defined on $[-1,1]$  has been study in  several papers (\cite{SS92,S14} and references therein) due to its importance to evaluate electron repulsion integrals \cite{S15}. Different Rys quadrature algorithms have been proposed in the literature \cite{K16,S14,S15}. Recently,   G. V.	Milovanović and N. Vasović, \cite{MV22}  considered the problem of construction of quadrature formulas of Gaussian type on $[-1, 1],$ concerning the following two-parametric weight function
$w(x; z) = e^{-zx^2}(1-x^2)^{\lambda -1/2}, \lambda >-1/2, $ which is a generalization of \textit{Rys quadrature} by  a transformation of quadratures on $[-1, 1]$ with $N$ nodes to ones on $(0, 1)$ with only $(N + 1)/2$ nodes. Such an approach provides a stable and very efficient numerical construction.\\

In this contribution we  explored the quadrature formula with respect to  the weight function $w(x; z) = e^{-zx}x^\alpha, \alpha >-1,$ supported on $[0,1]$. For this propose we use the \textit{Modified Chebyshev algorithm} using for this end the  \textit{shifted monic Jacobi} polynomials.  Here we make a computational analysis of our algorithm using for this objective the \texttt{QPQ} and \texttt{SQPQ} package of MATLAB provided by Walter Gautschi \cite{G18,G05}.  Finally,  we want to emphasize that the Truncated Laguerre weight appear in a natural way from the symmetrization process of  the weight $w(x;z)=|x|^{2\alpha+1} e^{-zx^2}$, which is  other type of generalization of \textit{Rys weights} \cite{DGM23}.

\section*{Acknowledgements}
The work of J. C. Garc\'ia-Ardila  has been supported by the Comunidad de Madrid multiannual agreement with the Universidad Rey Juan Carlos under the grant Proyectos I+D para Jóvenes Doctores, Ref. M2731, project NETA-MM. The work of F. Marcell\'an has been supported by FEDER/Ministerio de Ciencia e Innovación-Agencia Estatal de Investigación of Spain, grant PID2021-122154NB-I00, and the Madrid
Government (Comunidad de Madrid-Spain) under the Multiannual Agreement with UC3M in the line of Excellence of University Professors, grant EPUC3M23 in the context of the V PRICIT (Regional Program of Research and Technological Innovation).


\begin{thebibliography}{30}	
\bibitem{Chihara} T. S. Chihara, \textit{An Introduction to Orthogonal Polynomials}. In: Mathematics and its Applications Series, Vol. \textbf{13}. Gordon and Breach Science Publishers, New York-London-Paris, 1978.	
		
\bibitem{DM22} D. Dominici, F. Marcellán, \textit{Truncated Hermite polynomials}, J. Difference Equ. Appl. 2023. In press.	
		
\bibitem{DGM23} D. Dominici, J. C. García-Ardila, F. Marcellán, \textit{Symmetrization process and truncated orthogonal polynomials}. arXiv:2307.09581v2 [math.CA].
		
\bibitem{D76} M. Dupuis, J. Rys,  H. F. King. \emph{Evaluation of molecular integrals over Gaussian basis functions}. J. Chem. Phys. \textbf{65}(1), 111- 116. (1976).
		
\bibitem{GMM21} J. C. García-Ardila, F. Marcellán, M. E. Marriaga, \textit{Orthogonal Polynomials and Linear Functionals. An Algebraic Approach and Applications.} EMS Series of Lectures in Mathematics, EMS Press, Berlin, 2021.
		
\bibitem{G70} W. Gautschi, \textit{On the construction of Gaussian quadrature rules from modified moments.} Math. Comp. {\bf 24}, 245--260.	(1970)
		
\bibitem{Ga78} W. Gautschi, \textit{Questions of numerical conditions related to polynomials} in \emph{Symposium on Recent Advances in Numerical Analysis}, C. de Boor and G. H. Golub, Eds. Academic Press, New York, 1978,  pp. 45–72.
		
\bibitem{Ga82} W. Gautschi, \emph{An algorithmic implementation of the generalized Christoffel theorem}, in \emph{Numerical Integration} G. H\"{a}mmerlin Editor, International Series of Numerical Mathematics, Vol. \textbf{57}, Birkh\"{a}user Verlag,, Basel.1982
		
\bibitem{Ga02} W. Gautschi, \textit{The interplay between classical analysis and (Numerical) Linear Algebra- A tribute to Gene Golub}, ETNA \textbf{13}, 119--147. (2002).
		
\bibitem{G05} W. Gautschi, \textit{Orthogonal polynomials (in Matlab)}. J. Comput. Appl. Math. \textbf{178}, 215-234. (2005).
		
\bibitem{G83} W. Gautschi, \textit{Orthogonal polynomials, quadrature, and approximation: Computational methods and software (in Matlab)}, in \textit{Orthogonal Polynomials and Special Functions: Computation and Applications} (F. Marcellán and W. Van Assche, eds.), Lecture Notes in
Mathematics 1883, Springer, Berlin, 2006, pp. 1–77.
		
\bibitem{G16} W. Gautschi, \textit{Orthogonal Polynomials in Matlab: Exercises and solutions, Software Environment. Tools.} SIAM, Philadelphia, PA , 2016.
		
\bibitem{G18} W. Gautschi,  \textit{A Software Repository for Orthogonal Polynomials. Software, Environments and Tools.} {\bf  28.}  SIAM, Philadelphia, PA,  2018.
		
\bibitem{Go69} G. H. Golub,  J. H. Welsch, \textit{Calculation of Gauss quadrature rules}, Math. Comp., \textbf{23}, 221- 230. (1969).
		
\bibitem{GR80} I. S. Gradshteyn, I. M. Ryzhik, \textit{Table of Integrals, Series, and Products, Corrected and enlarged edition}, Academic Press, New York, First printing, 1980.		
		
\bibitem{MK05}	M. E. H. Ismail, E. Koelink,  \textit{Theory and applications of special functions. A volume dedicated to Mizan Rahman.} Papers from the Special Session of the American Mathematical Society Annual Meeting held in Baltimore, MD, January 15–18, 2003. Edited by Mourad E. H. Ismail and Erik Koelink Dev. Math., vol. 13 Springer, New York, 2005.
		
\bibitem{K16} H. F. King, \textit{ Strategies for evaluation of Rys roots and weights.} J. Phys. Chem. A {\bf 120}, 9348--9351. (2016).
		
\bibitem{Lyu} S. Lyu, Y. Chen,  \emph{Gaussian unitary ensembles with two jump discontinuities, {PDE}s, and the coupled Painlev\'{e} $II$  and $IV$ systems}. Stud. Appl. Math. \textbf{146} (1), 118--138. (2021).
		
\bibitem{M15} G. V. Milovanović, \textit{Construction and applications of Gaussian quadratures with nonclassical and exotic weight functions.} Stud. Univ. Babes-Bolyai Math. {\bf 60}, 211–233. (2015).
		
\bibitem{M18} G. V.  Milovanović, \textit{ An efficient computation of parameters in the Rys quadrature formula.} Bull. Cl. Sci. Math. Nat. Sci. Math. {\bf  43}, 39-64. (2018).
		
\bibitem{MMR94}  G. V. Milovanović, D. S. Mitrinović, Th. M. Rassias, \textit{Topics in polynomials: extremal problems, inequalities, zeros.} World Scientific Publishing Co., Inc., River Edge, NJ, 1994.
		
\bibitem{MV22} G. V. Milovanović, N. Vasović, \textit{Orthogonal polynomials and generalized Gauss-Rys quadrature formulae}. Kuwait J. Sci. {\bf 49} no.1, 17 pp. (2022),
		
\bibitem{OLB10} F. W. J. Olver, D. W. Lozier, R. F. Boisvert, C. W. Clark, editors. \textit{NIST Handbook of Mathematical Functions}. U.S. Department of Commerce, National Institute of Standards and Technology, Washington, DC. Cambridge University Press, Cambridge, 2010.		
		
\bibitem{RD83} J. Rys, M. Dupuis,  H. F. King, \textit{Computation of electron repulsion integrals using Rys quadrature method.} J. Comput. Chem. {\bf 4}, 154--157. (1983).
		
\bibitem{SD71} R. A. Sack,  A. F. Donovan,  \textit{An algorithm for Gaussian quadrature given modified moments.} Numer. Math. {\bf 18}, 465--478. (1971).
		
\bibitem{SS92} R. P. Sagar,  V. H. Smith, \textit{On the calculation of Rys polynomials and quadratures.} Int. J. Quant. Chem. {\bf 42}(4), 827--836. (1992).
		
\bibitem{S14} R. P. Schwenke, \emph{On the computation of high-order Rys quadrature weights and nodes}. Comput. Phys. Comm. {\bf 185}, 762–763. (2014).
		
\bibitem{S15} B. D. Shizgal,  \textit{A novel Rys quadrature algorithm for use in the calculation of electron repulsion integrals.} Comput. Theor. Chem. \textbf{1074}, 178–184. (2015).
		
\bibitem{W74} J. C. Wheeler, \textit{Modified moments and  Gaussian quadratures.} Rocky Mountain J. Math., {\textbf 4},   287–296. (1974).
		
\end{thebibliography}
\end{document}